\documentclass{article}

\usepackage{amsfonts}
\usepackage{amssymb}
\usepackage{amsmath}
\usepackage{amsthm}

\usepackage{clrscode}   
\usepackage{authblk}    

\usepackage[utf8]{inputenc} 

\usepackage{tikz}
\usetikzlibrary{patterns}

\newtheorem{theorem}{Theorem}[section]

\newtheorem{lemma}[theorem]{Lemma}

\begin{document}

\title{Flattening Karatsuba's recursion tree\\ into a single summation}
\date{February 2019\thanks{The first version of this paper was carefully reviewed by Aur\'elien Monteillet and Anthony Travers. Their comments have been taken into account in the current revision.}}
\author[$\dagger$]{Thomas Baruchel}
\affil[$\dagger$]{\small Éducation nationale, France\authorcr\texttt{baruchel@riseup.net}}
\maketitle

\begin{abstract}
\noindent The recursion tree resulting from Karatsuba's formula is built here by using an interleaved splitting scheme rather than the traditional left/right one. This allows an easier access to the nodes of the tree and some of them are initially flattened all at once into a single recursive formula. The whole tree is then flattened further into a convolution formula involving less elementary multiplications than the usual Cauchy product --- leading to iterative (rather than recursive) implementations of the algorithm. Unlike the traditional splitting scheme, the interleaved approach may also be applied to infinite power series, and corresponding formulas are also given.
\end{abstract}

\section{Introduction}
The fast multiplication algorithm discovered by Anatoly Karatsuba in 1960 (and published two years later) is known to be the oldest algorithm faster than the ``grade school'' method (when involved numbers or polynomials are large enough); while newer algorithms are still faster for sufficiently large numbers or polynomials, it is still widely used today for multiplicating medium-sized numbers or polynomials.

Due to its recursive divide-and-conquer approach, implementing this algorithm with no care about various issues (mostly related to storage of the temporary data) will lead to poor and often slow programs. Furthermore, the triple recursion involved by the algorithm, along with propagating changes in the computed data due to consecutive subtractions, makes implementing it in an iterative style more challenging.

The purpose of this paper is to deeply rewrite Karatsuba's formula in such a way that an iterative implementation would at first glance naturally arise. The last example of high-level pseudocode given in the paper thus relies on a single simple loop. The current paper will rather focus on identifying and writing down the new formula rather than on computational issues.

While keeping the very same number of elementary products, namely $n^{ \log _2 3}$ multiplications, where $n$ is the degree of the polynomials, the formulas that will be presented in section~\ref{full} do not necessarily keep also the number of additions and subtractions as low as in the state-of-the-art implementations of the algorithm --- deeper study of the question has not been made however.

%

\section{Karatsuba's recursion tree}\label{tree}
Let $A$ and $B$ be two polynomials in the same indeterminate~$x$; the divide-and-conquer paradigm to be used in the next sections requires splitting both~$A$ and~$B$ with the help of a third polynomial~$X$ (most likely a monomial)  in the same indeterminate~$x$; the exact purpose of~$X$ will be discussed in the current section as well as the following one.

While the ideas discussed here may be applied to any variant of Karatsuba's initial formula, we take the following one as a starting point and group all terms as factors around each of its three distinct branches:
\begin{equation}\label{karatsuba}
AB = 
    \underbrace{\left(X+1\right)A_0 B_0}_\textrm{branch 0}
  + \underbrace{X\left(X+1\right)A_1 B_1}_\textrm{branch 1}
  - \underbrace{X\left(A_1 - A_0\right)\left(B_1 - B_0\right)}_\textrm{branch 2}
\end{equation}
with $A=A_1X+A_0$ and $B=B_1X+B_0$. The formula is intended to be applied recursively until elementary products are encountered, and at each step the five terms $A_0$, $A_1$, $B_0$, $B_1$ and $X$ must be redefined according to the exact level of the recursion. Each node in the recursion tree will be labelled here according to a radix-3 labelling system --- by aggregating the arbitrary reference numbers, which are specified in~(\ref{karatsuba}), of the successive branches leading to it. Computing the product $AB$ now amounts to summing all contributions corresponding to the leaf nodes of a tree.

When a node is reached (by starting from the root node) without walking on any branch-2 nodes, we call it here a \emph{direct} node; it will otherwise be called \emph{indirect}. Obviously, direct nodes are reached by walking along paths whose label does not contain any digit~$2$. An indirect node whose parent is a direct node will be called \emph{primary indirect}. The following diagram illustrates that by showing direct and primary indirect nodes in a ternary recursion tree: hatched nodes are the primary indirect ones while all others are direct nodes (furthermore, shaded nodes are direct leaf nodes). Non-primary indirect nodes are discarded here, since they will be taken into account during the recursive process applied to each primary indirect node.

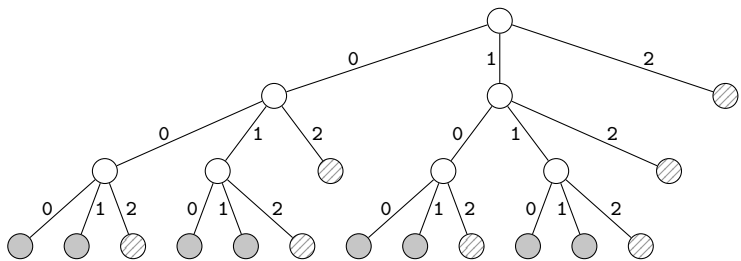
\begin{center}\begin{tikzpicture}
    \tikzstyle{mynode}=[circle, draw]
    \tikzstyle{leaf}=[fill=gray!45]
    \tikzstyle{indirect}=[pattern=north east lines, pattern color=gray!80]
    \tikzstyle{myarrow}=[]
    \node[mynode] (L1) at (0, 3) {};
    \node[mynode] (L2a) at (-3, 2) {};
    \node[mynode] (L2b) at (0, 2) {};
    \node[mynode,indirect] (L2c) at (3, 2) {};
    \node[mynode] (L3a) at (-5.25, 1) {};
    \node[mynode] (L3b) at (-3.75, 1) {};
    \node[mynode,indirect] (L3c) at (-2.25, 1) {};
    \node[mynode] (L3d) at (-0.75, 1) {};
    \node[mynode] (L3e) at (0.75, 1) {};
    \node[mynode,indirect] (L3f) at (2.25, 1) {};
    \node[mynode, leaf] (L4a) at (-6.375, 0) {};
    \node[mynode, leaf] (L4b) at (-5.625, 0) {};
    \node[mynode,indirect] (L4c) at (-4.875, 0) {};
    \node[mynode, leaf] (L4d) at (-4.125, 0) {};
    \node[mynode, leaf] (L4e) at (-3.375, 0) {};
    \node[mynode,indirect] (L4f) at (-2.625, 0) {};
    \node[mynode, leaf] (L4g) at (-1.875, 0) {};
    \node[mynode, leaf] (L4h) at (-1.125, 0) {};
    \node[mynode,indirect] (L4i) at (-0.375, 0) {};
    \node[mynode, leaf] (L4j) at (0.375, 0) {};
    \node[mynode, leaf] (L4k) at (1.125, 0) {};
    \node[mynode,indirect] (L4l) at (1.875, 0) {};
    \draw[myarrow] (L1) to node[midway,above=0pt, left=7pt]{\footnotesize\texttt{0}} (L2a);
    \draw[myarrow] (L1) to node[midway,above=0pt, left=-2.5pt]{\footnotesize\texttt{1}} (L2b);
    \draw[myarrow] (L1) to node[midway,above=0pt, right=8pt]{\footnotesize\texttt{2}} (L2c);
    \draw[myarrow] (L2a) to node[midway,above=0pt, left=4pt]{\footnotesize\texttt{0}} (L3a);
    \draw[myarrow] (L2a) to node[midway,above=0pt, right=-1pt]{\footnotesize\texttt{1}} (L3b);
    \draw[myarrow] (L2a) to node[midway,above=0pt, right=0pt]{\footnotesize\texttt{2}} (L3c);
    \draw[myarrow] (L2b) to node[midway,above=0pt, left=-0.5pt]{\footnotesize\texttt{0}} (L3d);
    \draw[myarrow] (L2b) to node[midway,above=0pt, left=-1pt]{\footnotesize\texttt{1}} (L3e);
    \draw[myarrow] (L2b) to node[midway,above=0pt, right=5pt]{\footnotesize\texttt{2}} (L3f);
    \draw[myarrow] (L3a) to node[midway,above=0pt, left=0mm]{\footnotesize\texttt{0}} (L4a);
    \draw[myarrow] (L3a) to node[midway,above=0pt, right=-2pt]{\footnotesize\texttt{1}} (L4b);
    \draw[myarrow] (L3a) to node[midway,above=0pt, right=-1pt]{\footnotesize\texttt{2}} (L4c);
    \draw[myarrow] (L3b) to node[midway,above=0pt, left=-1.5pt]{\footnotesize\texttt{0}} (L4d);
    \draw[myarrow] (L3b) to node[midway,above=0pt, left=-2.5pt]{\footnotesize\texttt{1}} (L4e);
    \draw[myarrow] (L3b) to node[midway,above=0pt, right=1pt]{\footnotesize\texttt{2}} (L4f);
    \draw[myarrow] (L3d) to node[midway,above=0pt, left=0mm]{\footnotesize\texttt{0}} (L4g);
    \draw[myarrow] (L3d) to node[midway,above=0pt, right=-2pt]{\footnotesize\texttt{1}} (L4h);
    \draw[myarrow] (L3d) to node[midway,above=0pt, right=-1pt]{\footnotesize\texttt{2}} (L4i);
    \draw[myarrow] (L3e) to node[midway,above=0pt, left=-1.5pt]{\footnotesize\texttt{0}} (L4j);
    \draw[myarrow] (L3e) to node[midway,above=0pt, left=-2.5pt]{\footnotesize\texttt{1}} (L4k);
    \draw[myarrow] (L3e) to node[midway,above=0pt, right=1pt]{\footnotesize\texttt{2}} (L4l);
\end{tikzpicture}\end{center}

Obviously, computing the whole product can also be seen as summing the contributions of all direct leaf nodes and primary indirect nodes (which are either shaded or hatched on the diagram above). This strategy will be deeper studied in the section~\ref{partial}.

\section{The art of splitting polynomials}\label{art}
Karatsuba's algorithm seems to be most of the time implemented or studied by using the same splitting scheme which involves taking apart terms of lower and higher degree (this may intuitively be seen as a left/right approach); it will be referred to here as the ``traditional splitting scheme''.

Several other splitting schemes are fully compliant with formula~(\ref{karatsuba}), namely any scheme taking apart groups of some power-of-$2$ sequential terms. The simplest one will be considered from now on: taking apart terms of even and odd rank. Such choice will have two main benefits: identifying the exact labelling number of a node is now easier and applying the algorithm to infinite power series will also be possible. This splitting scheme will be referred to here as the ``interleaved splitting scheme''. Applying it recursively is illustrated below:
\[
    \begin{array}{l}
        a_7x^7 + a_6x^6 + a_5x^5 + a_4x^4 + a_3x^3 + a_2x^2 + a_1x + a_0 \\[8pt]
       \quad=\left(a_7 x^6 + a_5x^4 + a_3 x^2 + a_1\right)x
        + \left(a_6x^6 + a_4x^4 + a_2x^2 + a_0 \right)\\[8pt]
        \quad=\left(\left(a_7 x^4 + a_3\right)x^2 + \left(a_5x^4 + a_1\right)\right)x
        + \left(\left(a_6x^4 + a_2\right)x^2 + \left(a_4x^4 + a_0 \right)\right)
    \end{array}
\]
New iterations of the interleaved splitting scheme lead to increasingly-sparse polynomials: each next term has initially a non-null coefficient, then each second term, then each fourth term, etc. If the term~$X$ in the formula~(\ref{karatsuba}) is initially some~$x$, then it will become~$x^2$ at the second iteration, then~$x^4$, etc.
\begin{lemma}\label{lemma}
    Let $K$ be a direct leaf node, built according to the interleaved splitting scheme, and reached by following the $(d_1d_2d_3\dots d_m)_3$-labelled path (where all $d_k$ are digits from $\{0,1\}$ since a direct node has no digit~$2$ in its path and with $m=\log_2 n$). The contribution of this node~$K$ to the whole summation is
\[
    \displaystyle\frac
    {1-x^{2^m}}
    {1-x}\,
    a_r b_r\, x^r
    \quad\textrm{with}\,\,r=(d_m d_{m-1} \dots d_3 d_2 d_1)_2
\]
which means that an initial radix-3 string is merely read later as a binary string with no further conversion (other than reverting the order of the digits).
\end{lemma}
\begin{proof}
    The left factor is the same for all direct leaf nodes; it does not rely on the path to the node~$K$ but only on the number~$m$ of iterations. Its purpose is to accumulate all iterated $X+1$ factors from the formula~(\ref{karatsuba}), namely
    \[(1+x)(1+x^2)(1+x^4)\dots=1+x+x^2+x^3+x^4+x^5+\dots x^{2^m-1}\,\textrm{.}\]
The fraction is a shorthand notation for this expression\footnote{The $1/(1-x)$ part is an usual generating function for~$1+x+x^2+x^3+x^4+\dots$ and multplicating it by~$(1-x^{2^m})$ allows to truncate the series to an arbitrary degree.}.

The $x^r$ part, which is the fourth factor in the formula from the lemma, comes by following the given path made of selected branches~0 and branches~1 in different levels of the tree. Accumulating the iterated~$X$ factors from the branch~1 in the formula~(\ref{karatsuba}) again is done as:
\[
    (x^1)^{d_1}\,
    (x^2)^{d_2}\,
    (x^4)^{d_3}\,
    \dots
    (x^{2^{m-1}})^{d_m}
\]
where the digit~$d_k$ is used for indicating whether the~$x^{2^{k-1}}$ factor is accumulated or not ---~which is only the case in a branch~1, allowing to use the digit as an exponent. Of course, the whole product is equal to~$x^r$.

Finally, the~$a_r$ and~$b_r$ factors are identified by induction: in the formula~(\ref{karatsuba}) $A$ and~$B$ are increasingly-sparse polynomials having both $a_k$ and $b_k$ as their constant term for some arbitrary direct node of the tree, where~$k$ also is the degree of all previously accumulated~$X$ factors. This is obviously true at the root level of the tree (with~$k=0$), and this property remains true at each new level of the tree whatever the selected branch is\footnote{This property can also be noticed in the illustration of the interleaved splitting scheme at the beginning of the current section for a polynomial of degree~$7$, formally applying this splitting scheme being deeply related to applying the formula~(\ref{karatsuba}) in regards to the degree of computed terms. In this example, the constant term of each sparse polynomial in parentheses at any level is $a_k$ when the external factor of this polynomial is $x^k$.}. Thus, the selected terms in~$A$ and in~$B$ are the coefficients of the terms of degree~$r$ in each polynomial.
\end{proof}


\section{Partially flattening the recursion tree}\label{partial}

In this section, we compute all direct leaf nodes at once while gathering separately all primary indirect nodes (thus taking care of all shaded or hatched nodes in the diagram from the section~\ref{tree}).

The Lemma~\ref{lemma} helps achieving the first part of this goal: gathering all $a_rb_rx^r$ terms is done by computing the \emph{termwise} product of $A$ and $B$ ---~which will be written down as $A(x)\odot B(x)$ --- for both polynomials and series. Furthermore, the following conventions will be used from now on:
\begin{itemize}
    \item explicit symbols for multiplication and convolution ($\times$ and $\ast$ respectively) for formulae which are applied recursively;
    \item an implicit notation if one of the factors is a polynomial ---~or the generating function of an integer sequence~--- containing only $0$ and $1$ coefficients.
\end{itemize}

Let $A$ and $B$ be two polynomials of degree~$n-1$ (for clarity, $n$ being a power of~$2$). Once Karatsuba's formula has been applied repeatedly $\log\left(n\right)$ times, matching coefficients in~$A$ and~$B$ are multiplied together for all direct leaf nodes of the recursion tree. We can group the $n$ elementary multiplications as:
\begin{equation}\label{part1}
\displaystyle\frac{1-x^n}{1-x}\left(A\left(x\right)\odot B\left(x\right)\right)
\end{equation}
where, from a computational point of view, the left factor does not involve a true multiplication but rather $n$ \textit{shift}/\textit{add} operations ---~which match the required additions in Karatsuba's formula.

We now have to gather all primary indirect nodes, as stated at the end of the section~\ref{tree} and at the beginning of the current section. There are $2^k$ of them on the level $k$ of the tree ($k=0$ corresponding to the first level under the root node) as can be noticed on the diagram from the section~\ref{tree}, and we will collect them in the so-called \textit{level-order}\footnote{A tree is traversed in level-order when each node on a level is visited before going to a lower level. Of course, we only care here about primary indirect nodes.}.

By using two nested summation symbols, we can ``iterate'' on all levels (there are $\log_2(n)$ of them) and then on all~$2^k$ primary indirect nodes in the $k^{\textrm{th}}$~level.

We now have to remember that $A$ and $B$ are sparse polynomials, splitting them according to the interleaved splitting scheme and then subtracting the lowest part from the highest is equivalent to the two consecutive steps:
\begin{itemize}
    \item multiplicating $A(x)$ and~$B(x)$ by $1-x^{2^{k-1}}$;
    \item keeping only terms of degree~$x^{0\cdot 2^k}$, $x^{1\cdot 2^{k}}$, $x^{2\cdot 2^{k}}$, $x^{3\cdot 2^{k}}$, etc.
\end{itemize}

The second task involves filtering a polynomial in order to keep specific terms and cancel all others; this is easily achieved by performing the termwise product of such polynomial with a relevant \emph{mask} being another polynomial whose coefficients are in~$\{0, 1\}$. Building this mask is achieved by using some convenient tools from the theory of generating functions: it has to be remembered first that $1/(1-x^{2^k})$ expands to $1+x^{2^k}+x^{2\cdot 2^k}+x^{3\cdot 2^k}+\dots$, which can be truncated and shifted as required with:
\[ \displaystyle\frac{1-x^n}{1-x^{2^k}} \, x^{2^{k-1}+j} \]
where $x^{2^{k-1}+j}$ gives the relevant ``offset'' for the $j^{\textrm{th}}$~primary indirect node at the~$k^{\textrm{th}}$ level of the tree. Again, the binary encoding of~$2^{k-1}+j$ is closely related to the radix-3 string labelling the path for reaching the parent of a given primary indirect node: we merely track all accumulated~$X$ factors by iterating with the formula~(\ref{karatsuba}) before reaching the primary indirect node; more precisely, we want to match the term of lowest degree still available in~$A$ and in~$B$ (since many terms have been discarded through the splitting process).

In order to follow very closely the formula~(\ref{karatsuba}), each termwise product should immediately be divided by $x^{2^{k-1}+j}$ (before multiplicating the two newly-built polynomials) because we actually want each polynomial to have a constant term while separately accumulating the required~$X$ factors; but the division can occur later since
\[
    x^{2^{k-1}+j} \left(
    \displaystyle\frac{P(x)}{x^{2^{k-1}+j}}\,\times\,\frac{Q(x)}{x^{2^{k-1}+j}}
\right)
= \frac{P(x)\times Q(x)}{x^{2^{k-1}+j}}
\]
where the three factors match those in the branch~2 in the formula~(\ref{karatsuba}), except for one thing: while all~$X$ factors have been correctly accumulated here, the~$X+1$ ones (from the branches~$0$ and~$1$ in the previous steps) are still missing ---~but there is exactly one of them at each level of the tree, as previously, and we can accumulate them with the same method than in the formula~(\ref{part1}).

Thus, gathering all primary indirect nodes finally gives:
\[
\begin{array}{l}
-
\displaystyle\sum_{k=1}^{\log_2\left(n\right)}
  \; \displaystyle\sum_{j=0}^{2^{k-1}-1}
  \; \displaystyle\frac{1-x^{2^{k-1}}}{\left(1-x\right)x^{2^{k-1}+j}} \\[18pt]
  \qquad  \quad   \displaystyle
  \left(
  \displaystyle\frac{1-x^n}{1-x^{2^k}} x^{2^{k-1}+j} \odot
\left(1-x^{2^{k-1}}\right)A\left(x\right)
\right) \\[18pt]
\qquad \times
\left(
  \displaystyle\frac{1-x^n}{1-x^{2^k}} x^{2^{k-1}+j} \odot
\left(1-x^{2^{k-1}}\right)B\left(x\right)
\right)
\end{array}
\]
where $k$ is the level of each considered row of nodes, and $j$ the rank of each \emph{direct} subtracting node on the level~$k$.

Since $2^{k-1}+j$ merely iterates over $1,2,3,\dots,n-1$, it is easy to use a single summation for directly iterating over all the considered subtrees:
\begin{equation}\label{part2}
\begin{array}{l}
-
\displaystyle\sum_{m=1}^{n-1}
\displaystyle\frac{1-x^{2^{\lfloor\log_2(m)\rfloor}}}{\left(1-x\right)x^{m}} \\[18pt]
  \qquad  \quad   \displaystyle
  \left(
      \displaystyle\frac{\left(1-x^n\right) x^m}{1-x^{2^{\lfloor\log_2(m)\rfloor+1}}} \odot
      \left(1-x^{2^{\lfloor\log_2(m)\rfloor}}\right)A\left(x\right)
\right) \\[18pt]
\qquad \times
\left(
    \displaystyle\frac{\left(1-x^n\right) x^m}{1-x^{2^{\lfloor\log_2(m)\rfloor+1}}} \odot
    \left(1-x^{2^{\lfloor\log_2(m)\rfloor}}\right)B\left(x\right)
\right)
\end{array}
\end{equation}

Summing both parts (\ref{part1}) and (\ref{part2}) results in $A\left(x\right)\times B\left(x\right)$ and is more or less equivalent to Karatsuba's algorithm from a computational point of view ---~as long as some variable substitution is done before each recursive call in order to map sparse polynomials to new polynomials of smaller degree.

Using the previously described interleaved splitting scheme now allows to apply Karatsuba's recursive formula to infinite power series (which is not the case with the traditional splitting scheme).

In the formulas~(\ref{part1}) and~(\ref{part2}), all $1-x^n$ numerators in the generating functions are intended to truncate periodical sequences of unitary and null coefficients to the required length. Extending these formulas to infinite power series is then easily achieved by removing such numerators:
\begin{equation}\label{series}
\begin{array}{l}
    f(x)\ast g(x)=
    \displaystyle\frac{f\left(x\right)\odot g\left(x\right)}{1-x}
-
  \displaystyle\sum_{m=1}^\infty
  \; \displaystyle\frac{1-x^{2^{\lfloor\log_2 m\rfloor}}}{\left(1-x\right)x^m} \\[18pt]
  \qquad\qquad\qquad\qquad\qquad\begin{array}{ll}
      &\displaystyle
  \left(
      \displaystyle\frac{x^m}{1-x^{2^{\lfloor\log_2 m\rfloor+1}}}  \odot
\left(1-x^{2^{\lfloor\log_2 m\rfloor}}\right)f\left(x\right)
\right)\\[18pt]
\ast
      &\displaystyle
\left(
      \displaystyle\frac{x^m}{1-x^{2^{\lfloor\log_2 m\rfloor+1}}}  \odot
\left(1-x^{2^{\lfloor\log_2 m\rfloor}}\right)g\left(x\right)
\right)
\end{array}
\end{array}
\end{equation}

As a separate question, we may wonder for which values the index of summation~$m$ in the formulas~(\ref{part2}) and~(\ref{series}) will contribute to the computation of a term of degree~$d$ in the final result. Otherwise said, we want to gather the primary indirect nodes involved in a given resulting term ---~these primary indirect nodes are labelled according to the previously specified enumeration $1,2,3,\dots,n-1$ by following the level-order (on the diagram of the section~\ref{tree}, the seven hatched nodes will be labelled~$1,2,3,\dots,7$). The set~$S_d$ of all such indices~$m$ is:
\begin{equation}\label{degree}
    S_d = \left\{  m \,\Big| \, 1\leqslant m\leqslant d,  \left(\left(d-m\right) \,\textrm{mod}\,\, 2^{\lfloor\log_2 m\rfloor+1}\right) < 2^{\lfloor\log_2 m\rfloor}\right\}
\end{equation}
which comes directly from the formula~(\ref{part2}): the $m^{\textrm{th}}$ primary indirect node have some~$x^m P(x)$ contribution ---~with some $P(x)=c_0 + c_1 x^{2^{k+1}} + c_2 x^{2\cdot 2^{k+1}} + \dots$ (at the $k^{\textrm{th}}$~level of the tree, and of course $k=\lfloor\log_2 m\rfloor$). Because of all~$X+1$ accumulated factors, each computed coefficient will also be shifted $2^k$ times ``to the right''. The definition~(\ref{degree}) gathers all primary indirect nodes such that one term in the involved sparse polynomial has a degree ``close'' enough to~$d$.

The cardinality $\vert S_d\vert$ of such sets of indices is empirically found to be the sequence \texttt{A268289} in the \textit{On-Line Encyclopedia of Integer Sequences}, namely the cumulated differences between the number of digits 1 and the number of digits~0 in the binary expansions of integers up to~$d$\hspace{0.15em}. Another explicit expression for $\vert S_d\vert$ resorting to the $\tau$ Takagi function can also be given:
\begin{equation}
\vert S_d\vert = 
        {\normalfont\texttt{A268289}}_d
        =
        d - 2^k \tau\left(\displaystyle\frac{d+1} {2^k }   - 1  \right)
\end{equation}
with $d$ some non-negative\footnote{By convention ${\normalfont\texttt{A268289}}_0=0$, the formal definition of the sequence being slightly different than the plain english one above which would imply the wrong statement ${\normalfont\texttt{A268289}}_0=-1$.} integer and $k=\lfloor\log_2(d)\rfloor$.

While inserting an extended proof of the previous identity would be far beyond the scope of the current paper, a quick hint will help building such a proof: when terms of the three sequences for all indices up to~$2^s-1$ are known, we build the following terms up to the index~$2^{s+1}-1$ with the help of the same building rule for the three sequences:
\[
    u_{n+2^s}
    = u_n
        + \left(n+1\right) \left(\lfloor\log_2(n)\rfloor-s+2\right) + 2^s - 2^{\lfloor\log_2(n)\rfloor+1}
\]
and since the three sequences share the same initial terms, we finally prove that they are identical.
%

\section{Fully flattening the recursion tree}\label{full}
Having described in the previous section how to handle two branches of the tree at once by using the termwise multiplication formula, we now go one step further. The key ideas from the previous section obviously are:
\begin{itemize}
    \item handling several leaf nodes as a whole by using termwise products of polynomials;
    \item performing all required operations (additions, subtractions, shifts) by accumulating factors finally expanding as polynomials with unitary and null coefficients;
    \item keeping or cancelling coefficients of polynomials by performing the termwise product of the latter with relevant masks.
\end{itemize}

We will extend the application of these ideas in order to handle all the $n^{\log_2 3}$~leaf nodes by computing only $n$~termwise products. As a starting point, we consider again the interleaved splitting scheme which is illustrated in Section~\ref{art} by the following example:
\[
    \begin{array}{l}
        a_7x^7 + a_6x^6 + a_5x^5 + a_4x^4 + a_3x^3 + a_2x^2 + a_1x + a_0 \\[8pt]
       \quad=\left(a_7 x^6 + a_5x^4 + a_3 x^2 + a_1\right)x
        + \left(a_6x^6 + a_4x^4 + a_2x^2 + a_0 \right)\\[8pt]
        \quad=\left(\left(a_7 x^4 + a_3\right)x^2 + \left(a_5x^4 + a_1\right)\right)x
        + \left(\left(a_6x^4 + a_2\right)x^2 + \left(a_4x^4 + a_0 \right)\right)
    \end{array}
\]
The last line of the example shows four sparse polynomials which would occur in four different paths, namely $(11)_3$, $(10)_3$, $(01)_3$ and $(00)_3$ (when reading these paths in the reversed order as binary encoded strings, we can identify the index of the constant term for all these polynomials). Instead of exploring again the branches~0 and~1, we now explore the branch~2 from each of these four nodes ---~reaching the lowest level of the tree with four elementary terms:
\[
        \left(a_7 - a_3\right) , \left(a_5 - a_1\right),
        \left(a_6 - a_2\right) \,\textrm{and}\, \left(a_4 - a_0 \right)
\]
which have to be multiplicated with corresponding~$b$ coefficients. We want to compute the four multiplications as~$\left(1-x^4\right)A(x)\odot\left(1-x^4\right)B(x)$. Extraneous coefficients are easily cancelled with a mask. While the mask should be applied before the termwise product from a computational point of view, we focus rather on building the most concise formula and the use of the mask will be postponed.

It has to be noticed that the~$1-x^{2^{k-1}}$ factor is the same for all nodes from the $k^{\textrm{th}}$ level; of course such factors may also be accumulated when walking on several branches~$2$ on a given path. Since all nodes must be visited by the recursion process, all possible selections of such factors have to be considered. The ternary initial tree built from the formula~(\ref{karatsuba}) now becomes a more classical binary tree: at each level we can choose between either accumulating a new~$1+X$ factor (to be used \emph{after} the termwise product) or accumulating a new~$1-X$ factor (to be used \emph{before} the termwise product). Indeed, walking on a branch~$0$ or~$1$, implies accumulating the relevant~$1+X$ factor according to the initial identity~(\ref{karatsuba}), while walking on the branch~$2$ implies subtracting coefficients to other ones wich is performed here by accumulating the relevant~$1-X$ factor. The~$1-X$ factors have an arithmetical purpose and must be applied on the actual numerical values \emph{before} the termwise multiplication (which is going to definitely discard some coefficients), while the~$1+X$ factors have a shifting/adding purpose and must be applied \emph{after} the termwise purpose when the unwanted coefficients have been discarded.

Accumulating various~$1+X$ factors actually plays two different roles in the computation: one has been already described (shifting and adding some terms), we focus now on the other one. The mask to be applied for discarding extraneous and useless subtracted coefficients happens to be the very same polynomial made of accumulated~$1+X$ factors as long as we shift and truncate it accordingly.

This is proved as follows; when computing the product
\[
    \left(1+x\right) \left(1+x^2\right) \left(1+x^4\right) \left(1+x^8\right) \dots
\]
and considering a specific factor~$1+x^{2^k}$, we have to remember that the latter gives some control on each block of $2^{k+1}$~consecutive coefficients: either the left pattern of $2^k$~coefficients is duplicated on the right part or not; for instance, removing the single factor~$1+x^2$ would lead to the sequence of coefficients~$1,1,0,0,1,1,0,0,1,1,0,0,\dots$ When ``tracking'' the coefficients during the recursion process (according to the interleaved splitting scheme), the $1-x^{2^k}$ factor means subtracting the left part of such a block from the right part, and of course we must cancel half of the coefficients in each block.

We can now write down a formula for the whole tree by choosing some way to iterate over all selections of~$1-X$ factors or over all selections of~$1+X$ factors; this could be done, for instance, by iterating over all subsets of~$\{1, 2, \dots, \log_2 n\}$. By noticing that the mask is being iterated over all divisors of an appropriate polynomial, we can utilise some extra conventions and notations to arrive at an elegant formulation:
\begin{equation}\label{formula}
    \begin{array}{r}
        A\times B \, = \!\!\!\!\!\!\!{\displaystyle\sum_{
    \begin{array}{c}
        \text{\footnotesize$f\in\mathbb{Z}[X],$}
        \\
        \text{\footnotesize$f \,\big |\, \sum_{k=0}^{n-1}X^k$}
     \end{array}}}\!\!\!\!\!\!
     \displaystyle
     f\left(f \,\dot{f}^\textrm{\tiny\hspace{1pt}W} \odot \dot{f}A\odot\dot{f}B\right)
 \end{array}
\end{equation}
\noindent with $\dot{f}$ selecting all unselected $1+X^{2^k}$ factors in $f$ and negating ---~for each one~--- the coefficient of their non-constant term, and $\dot{f}^\textrm{\tiny\hspace{1pt}W}$ the leading term (including the coefficient) of $\dot{f}$. The superscript character W stands for ``weight''\footnote{This symbol is compact but not very common; it can be found however in an article by Shigeru Kuroda, \emph{Shestakov-Umirbaev reductions and Nagata's conjecture on a polynomial automorphism} (2007).}.


The divisors of the $\sum_{k=0}^{n-1}X^k$ polynomial are of course all possible
selections of factors in~$(1+X)(1+X^2)(1+X^4)\dots$ The
$\dot{f}^\textrm{\tiny\hspace{1pt}W}$ term is used for shifting the mask at the
beginning of the meaningful subtracted terms, and the coefficient of this term
is either~$-1$ or~$1$ according to the number of subtractions (respectively odd or even) in the current
path.

Of course the expected~$n^{\log _2 3}$ elementary multiplications are embedded in the previous formula, since the whole idea of the current section was to track them and to \textit{shift/add/subtract} them according to the initial~(\ref{karatsuba}) formula; no supplementary multiplication was added anywhere.

An alternate version of the formula comes from the fact that $\dot{f}^\textrm{\tiny\hspace{1pt}W}$ is either some~$-x^k$ or~$x^k$ with~$k$ being the sum of powers of~$2$ involved in the product of~$1-x^{2^j}$ when building the~$\dot{f}$ polynomial. We can thus iterate over~$1, 2, \dots, n-1$ as the degree of $\dot{f}^\textrm{\tiny\hspace{1pt}W}$ and build other terms from it:
\begin{equation}\label{formula2}
    A\times B\, = \,
    \sum_{k=0}^{n-1}
    \bar{f}_k\left(\sigma_k \bar{f}_k X^k \odot \dot{f}_k A\odot\dot{f}_k B\right)
\end{equation}
\noindent with $\dot{f}_k=(1-X)^{d_0}(1-X^2)^{d_1}(1-X^4)^{d_2}\dots$ by referring to the binary digits of $k=(\dots d_2 d_1 d_0)_2$, with also~$\bar{f}_k=(1+X)^{1-d_0}(1+X^2)^{1-d_1}(1+X^4)^{1-d_2}\dots$ where the bar is intended to show that the product is truncated to the same ``format'' than~$\dot{f}_k$, made of $\log_2 n$ different factors, despite the infinite number of leading zeros in the binary encoding of~$k$, and with~$\sigma_k={\normalfont\texttt{A106400}}_k$ in the \textit{On-Line Encyclopedia of Integer Sequences} ($\sigma_k=-1$ if the binary weight of~$k$ is odd and $\sigma_k=1$ otherwise). Several expressions for~$\sigma_k$ are published on the page of the sequence~\texttt{A106400} ---~one involving an hypergeometric ${}_2F_1$ function.

Formula~(\ref{formula2}) can now easily be adapted to infinite power series as:
\begin{equation}\label{formula3}
    a(x)\ast b(x)\,=\,
    \sum_{k=0}^{\infty}
    f_k(x)\left(\sigma_k \,f_k(x)\, x^k \odot \dot{f}_k(x) \,a(x)\odot\dot{f}_k(x) \,b(x)\right)
\end{equation}
\noindent by merely removing the bar from~$\bar{f}_k$, now defining~$f_k$ as an infinite product.

\section{Implementing the new formula}\label{implementation}

Implementing the formula~(\ref{formula}) with no care about the true purpose of each part does not lead to a very efficient code: some elementary multiplications would be done though they are going to be cancelled soon after, many operations involving null coefficients in very sparse polynomials could be avoided, etc. But a first attempt can be given as a proof of concept and we give below two pieces of pseudocode intended to be used with any computer algebra system handling the polynomial type; they do not focus on low-level implementation issues (how more or less sparse polynomials are internally represented in order to give the most efficient access to their coefficients). Termwise multiplication of polynomials should of course be already implemented.

Two polynomials~$f$ and~$f'$ are used below for accumating \textit{shift/add} and \textit{shift/subtract} operations according to what was previously discussed.

The following code shows how subtracting and masking factors are accumulated while iterating on the branches of the tree; the key idea is to use the $1+X$ factors for selecting subtracted terms as well as for propagating them and the $1-X$ corresponding factors for performing the subtractions:
{\small
\begin{codebox}
    \Procname{$\proc{Multiply}(A, B)$}
    \li $d \gets \lceil\log_2(1+\max(\deg A, \deg B))\rceil$
    \li $n \gets 2^d$
    \li $s \gets 0$
    \li \For $k\gets 0$ \To $n-1$
    \li \Do
            $f\gets 1$
    \li     $f'\gets 1$
    \li     \For $j\gets 0$ \To $d-1$
    \li     \Do
                \If $k \mathrel{\&} 2^j \neq 0$ \>\>\>\>\>\>\>\Comment test if bit $j$ of $k$ is set
    \li         \Then $f'\gets (1-X^{2^j})\,f'$ \>\>\>\>\> \Comment subtracting factor
    \li         \Else $f\gets (1+X^{2^j})\,f$ \>\>\>\>\> \Comment termwise mask
                \End
            \End
    \li     $\id{degree}\gets \deg f'$
                     \>\>\>\>\>\>\>\>\>\Comment degree of $f'$
    \li     $\id{lt}\gets \proc{coeff}(f', X, \id{degree})\,X^{\id{degree}}$
                     \>\>\>\>\>\>\>\>\>\Comment leading term in $f'$
                     \li     $s\gets s + f \, (f\id{lt}\odot\, f'A\odot f'B)$
    \End
    \li \Return s
\end{codebox}
}
\noindent where it can be seen that the variable~$f$ has two distinct purposes: the algebraic one and also a tracking purpose for identifying how many terms have to be kept in the termwise multiplication, as explained in the previous section.

The previous pseudocode however loses all benefits of traditional implementations of Karatsuba's algorithm because the same subtractions are computed for distinct values of~$k$. Fortunately iterating over the binary expansions of~$k$ by using the reflected binary code (Gray code) instead of the standard radix-$2$ labelling system preserves the required number of subtractions.

Since iterating over such binary expansions (the Gray code) implies flipping a single bit between two consecutive integers, we are now visiting the nodes in a new unintuitive order in such a way that the accumulated factors can be reused from a node to another one by performing two single changes only: dividing the current product by some $1-X^k$ or by some~$1+X^k$ (according to the exact location of the flipped bit) in order to ``cancel'' the branch being left, as well as multiplicating it by some $1+X^k$ or by some~$1-X^k$ (according to the exact location of the flipped bit) in order to ``enable'' the new visited branch.

Of course, polynomial divisions are exact here since we merely divide some $(1+X)(1+X^2)(1+X^4)\dots$ polynomial by one of its previously ``enabled'' factors for cancelling it ---~the same for another $(1-X)(1-X^2)(1-X^4)\dots$ polynomial. Actually, both polynomial are always kept synchronized in the following way: if some $(1\pm X^k)$ factor is ``enabled'' in one polynomial, it is disabled in the other one (because both branches are not visited together at the same level).

The following version, though not optimized by itself from an implementation point of view (because it still relies on high-level polynomial types), gives the prototype of a more optimized iterative version of Karatsuba's algorithm:
{\small
\begin{codebox}
    \Procname{$\proc{Multiply2}(A, B)$}
    \li $d \gets \lceil\log_2(1+\max(\deg A, \deg B))\rceil$
    \li $n \gets 2^d$
    \li $f \gets (1-X^n)/(1-X)$
    \>\>\>\>\>\>\>\>\>\> \Comment initial mask (all bits set)
    \li $f' \gets 1$
    \li $g \gets 0$
          \>\>\>\>\>\>\>\>\>\> \Comment Gray-code counterpart of $k$
    \li $s \gets f\,(f\odot A\odot B)$
          \>\>\>\>\>\>\>\>\>\> \Comment case $g=k=0$
    \li \For $k\gets 1$ \To $n-1$
    \li \Do
    $j\gets \lfloor\log_2(k\;\textrm{xor}\;k-1)\rfloor$
          \>\>\>\>\>\>\>\> \Comment least significant set bit in $k$
    \li   
    \If $g \mathrel{\&} 2^j = 0$ \>\>\>\>\>\>\>\>\Comment test if bit $j$ has to be set in $g$
    \li         \Then
                      $f'\gets (1-X^{2^j})\,f'$ \>\>\>\>\>\> \Comment subtracting factor
    \li               $f\gets f/(1+X^{2^j})$ \>\>\>\>\>\> \Comment termwise mask
    \li         \Else
                      $f'\gets f'/(1-X^{2^j})$ \>\>\>\>\>\> \Comment subtracting factor
    \li               $f\gets (1+X^{2^j})\,f$ \>\>\>\>\>\> \Comment termwise mask
                \End
    \li     $g\gets g\;\textrm{xor}\;2^j$
    \>\>\>\>\>\>\>\>\Comment update $g$ (Gray-code of $k$)
    \li     $\id{degree}\gets \deg f'$
                     \>\>\>\>\>\>\>\>\Comment degree of $f'$
    \li     $\id{lt}\gets \proc{coeff}(f', X, \id{degree})\,X^{\id{degree}}$
                     \>\>\>\>\>\>\>\>\Comment leading term in $f'$
                     \li     $s\gets s + f \, (f\id{lt}\odot\, f'A\odot f'B)$
    \End
    \li \Return s
\end{codebox}
}
A lower-level implementation of this pseudocode should avoid actually storing the~$f'$ polynomial in a separate buffer and computing both $f'A$ and $f'B$ products ---~the idea being rather to directly store $f'A$ and $f'B$, and merely update them at each step of the loop.

Furthermore, efficiently implementing the previous pseudocode should take care of the subtracting and adding steps: since polynomials become very sparse for some values of~$k$, very few terms should be manipulated at these points. Two main directions should be explored for that purpose: using linked lists for representing polynomials or tracking the remaining non-null coefficients by using an elaborated system of strides\footnote{This is one of the most important concepts behind the famous \emph{Numpy} module for Python; strides allow to build views on parts of an existing array without actually copying it.}. Elementary multiplications should of course be aware of the mask to be applied in order to avoid useless computation.

\section{Computing arbitrary coefficients}

Adapting the formula~(\ref{formula3}) to a more or less Cauchy-like one is actually achievable ---~though very inefficient from a computational point of view since the same multiplications will be performed again and again rather than propagated through the \emph{shift/add} process. It may however have some interest in further theoretical investigations.

For that purpose, we do not need the $\odot$ termwise operator any longer since we have no interest in computing several coefficients at once, but we now need the $\&$ bitwise multiplication operator (the bitwise ``and'' operator) since the new formula will highly rely on testing whether such or such $1\pm X$ factor is selected or not.

Let $g(x)=a_0+a_1x+a_2x^2+\dots$ and $h(x)=b_0+b_1x+b_2x^2+\dots$, then
\[
    \begin{array}{l}
    \displaystyle g(x)\ast h(x)\,=\,
    \sum_{m=0}^\infty
    x^m
      \sum_{k=0}^m
      \sigma_k \\[12pt]
      \displaystyle\qquad\qquad\qquad
      \sum_{j=k}^m
      \tau_k(m,j) \,\tau_k(j,k)
        \left(\sum_{t=j-k}^j \upsilon_k(j,t)\, \sigma_t a_t \right)
        \left(\sum_{t=j-k}^j \upsilon_k(j,t)\, \sigma_t b_t \right)
    \end{array}
\]
with $\tau_k(m,j)=[k\mathbin{\&}(m-j)=0]$ and $\upsilon_k(j,t)=[k\mathbin{\&}(j-t)=j-t]$, both defined by using Iverson bracket, and $\sigma_k$ as previously defined in the formulas~(\ref{formula2}) and~(\ref{formula3}).

The $\tau_k$ function is used for testing whether a coefficient from a given degree will actually be shifted and added as a contribution to another given degree; the $\upsilon_k$ function is used for testing whether a coefficient of a given degree will be shifted and subtracted when accumulating all $1-X$ factors.

After having noticed that $\tau_k$ is related to the sequence \texttt{A047999} and that $\upsilon_k$ is related to the sequence \texttt{A106344}, we decide to clean up the previous formula by relying on \texttt{A047999} only, namely on the Sierpi\'nski triangle, finally getting:
\[
    \begin{array}{l}
    \displaystyle g(x)\ast h(x)\,=\,
    \sum_{m=0}^\infty
    x^m
      \sum_{k=0}^m
      \sigma_k \\[12pt]
      \displaystyle\qquad\qquad\qquad\qquad
      \sum_{j=k}^m
      T(k+m-j,k) \, T(j, k) \\[12pt]
        \displaystyle\qquad\qquad\qquad\qquad\qquad
        \left(\sum_{t=j-k}^j T(k, j-t)\, \sigma_t a_t \right)
        \left(\sum_{t=j-k}^j T(k, j-t)\, \sigma_t b_t \right)
    \end{array}
\]
with $T(n,k)$ being \texttt{A047999} --- defined as $T(n,k)=[k\mathbin{\&}(n-k)=0]$. While the resulting nested summations may look rather heavy, a quick glance at the famous graphical representation of Sierpi\'nski triangle will show that most of the involved terms are null.

\section{Conclusion}
While the section~\ref{implementation} highly relies on pseudocode parts, identifying and publishing the formulas~(\ref{formula}) to (\ref{formula3}) as new theoretical convolution formulas actually was the true purpose of the current paper. While the number of elementary multiplications remains identical to the expected~$n^{\log _2 3}$ one, actually applying the formula probably involves more additions and subtractions than what would be the case by following the conventional recursive approach --- their exact amount being however not investigated here.

The final formula, in the previous section, shows that tracking individual coefficients through the whole recursion tree is achievable; as a sophisticated convolution involving the Sierpi\'nski triangle, it may be seen as a starting point for further investigations focused towards combinatorics.

\vspace{32pt}
\noindent{\small \textbf{Conflict of Interest:} The authors declare that they have no conflict of interest.}

\vspace{8pt}
\noindent{\small The current article is accessible on \texttt{http://export.arxiv.org/pdf/1902.08982}\hspace{2pt}.}

\end{document}